 \documentclass[10pt]{article}
\usepackage{cite}
\usepackage{amsmath,amssymb,amsfonts}
\usepackage{algorithmic}
\usepackage{graphicx}
\usepackage{textcomp}
\usepackage{wrapfig}
\input{mysymbol.sty}
\usepackage[utf8]{inputenc}   
\usepackage[dvipsnames]{xcolor}
\usepackage{amsthm}

\usepackage{tikz}
\usetikzlibrary{shapes.geometric, arrows}
\tikzstyle{ad} = [rectangle, rounded corners, minimum width=3cm, minimum height=1cm,text centered, draw=black, fill=blue!30]
\tikzstyle{ag} = [rectangle, rounded corners, minimum width=2cm, minimum height=1cm,text centered, draw=black, fill=green!30]
\tikzstyle{cd} = [rectangle, minimum width=3cm, minimum height=1cm, text centered, draw=black, fill=orange!30]

\renewcommand{\baselinestretch}{0.98}\selectfont

\newtheorem{theorem}{Theorem}

\newtheorem{lemma}{Lemma}
\newtheorem{proposition}{Proposition}
\newtheorem{example}[theorem]{Example}
\newtheorem{definition}{Definition}
\newtheorem{remark}{Remark}
 
\begin{document}
\title{Maximizing Social Welfare and Agreement via Information Design in Linear-Quadratic-Gaussian Games}
\author{Furkan Sezer, Hossein Khazaei,  and Ceyhun Eksin,
\thanks{This work was supported by NSF CCF-2008855. Furkan Sezer and Ceyhun Eksin are with the Department of Industrial and Systems Engineering, Texas A\&M University, College Station, TX 77843 USA (e-mails: furkan.sezer@tamu.edu, eksinc@tamu.edu). Hossein Khazaei is with Department of Electrical and Computer Engineering, Stony Brook University 100 Nicolls Rd, Stony Brook, NY 11794 (e-mail: hossein.khazaei@stonybrook.edu) }}

\maketitle

\begin{abstract}
Information design in an incomplete information game involves a designer that aims to influence players' actions through signals generated from a designed probability distribution to optimize its objective function. For quadratic design objective functions, if the players have quadratic payoffs that depend on the players' actions and an unknown payoff-relevant state, and signals on the state that follow a Gaussian distribution conditional on the state realization, the information design problem is a semi-definite program (SDP) \cite{Ui2020}. In this note, we seek to characterize the optimal information design analytically by leveraging the SDP formulation, when the design objective is to maximize social welfare or the agreement among players' action.
We show that full information disclosure maximizes social welfare when there is a common payoff state, the payoff dependencies among players' actions are homogeneous, or when the signals are public. When the objective is to maximize the agreement among players' actions, not revealing any information is optimal. When the objective is a weighted combination of social welfare and agreement terms, we establish a threshold weight below which full information disclosure is optimal under public signals for games with homogeneous payoffs. Numerical results corroborate the analytical results, and identify partial information disclosure structures that are optimal. 
%
\end{abstract}

\section{Introduction}
Information design refers to the determination of the information fidelity of the signals given to the players in an incomplete information game, so that the induced actions of players maximize a system level objective \cite{Bergemann2019}. In an incomplete information game, players compete to maximize their individual payoffs that depend on the action of the player, other players' actions and an unknown state. Incomplete information games are used to model power allocation of users in wireless networks with unknown channel gains \cite{bacci2015game}, traffic flow in communication or transportation networks \cite{wu2021value}, consumer behavior in smart grids \cite{eksin2015demand}, coordination of autonomous teams \cite{eksin2014bayesian}, and currency attacks of investors \cite{morris2002social}. 
In such scenarios, we envision the existence of an information designer that can provide the ``best'' information about the payoff-relevant states to the players according to its objective (see Fig. \ref{inf_design_diagram}). 
As per the above examples, the information designer can represent an entity such as a system designer overseeing the spectrum allocation/traffic, a market-maker, an independent system operator in the power grid, or the federal reserve \cite{tavafoghi2017informational,sekar2019uncertainty,wu2021value}. Recently, information design is used in robust sensor design \cite{sayin_basar_2021}, and modeling deception/privacy \cite{sayin_basar_2022}.
In addition, information design framework is used to identify key/central players in social networks with respect to the goals of the system designer \cite{candogan2020social}, to maximize the utility of the insured agents in a competitive insurance market \cite{garcia2021information}, and to design public health warning policies against recurrent risks \cite{alizamir2020healthcare}. In an effort to obtain analytical solutions, and thus insights about the effects of the information design on the system, current approaches make structural assumptions about the state/action space, the system designer's objective, and the game payoffs---see  \cite{candogan2020information} for a detailed discussion. 

In this note, we study information design linear-quadratic-Gaussian (LQG) games.
In an LQG game, players have quadratic payoff functions, and the state and the signals come from a Gaussian distribution\begin{footnote}{We note that information design in LQG games is categorically different than the well-known LQG control which considers control of linear dynamical systems with Gaussian noise given quadratic control objectives.}\end{footnote}.
Under certain assumptions, the optimal strategy in LQG games defined by the Bayesian Nash equilibrium (BNE) is unique and linear in the signals received  \cite{radner1962}. The linearity of BNE strategies allow the information design problem to be a semi-definite program (SDP) when the information designer's objective is a quadratic function of the players' actions and the payoff-relevant states \cite{Ui2020} (Section \ref{sec_model}). 
Building on the SDP formulation of the information design problem, we analyze optimal information structures when the system level objective is to maximize social welfare (Section \ref{sec_welfare}), maximize agreement among players' actions (Section \ref{sec_conformism}), or a weighted combination of these two objectives (Section \ref{sec:mult}). 

%

An information structure comprises signal transmission rules and the probability distribution from which signals are generated. Signals transmitted to players convey information about payoff relevant states. 
The information structure is {\it public} when all players receive a common signal. Otherwise, when the players receive individual signals, the signal structure is {\it private}. Another distinction is based on the fidelity of information carried by the signals. A signal can carry {\it no}, {\it partial}, or {\it full information}. No information disclosure does not improve the prior information of the players about the payoff relevant state, while signals reveal the payoff relevant state under full information disclosure. A partial information disclosure is when the signals carry some information, but do not fully reveal the payoff relevant state to the players. 

\noindent {\bf Contributions:} In this paper, we provide analytical and computational insights about the value of information and optimal information structures by focusing on particular objectives for the designer (social welfare and agreement). 
{\bf 1)} Given the social welfare design objective, we show that full information disclosure is optimal if there is a common payoff state (Proposition \ref{prop_comm}), when the dependency of payoffs on others' actions is homogeneous (Theorem \ref{thm:sub-super}), or if we only consider the set of public signals (Proposition \ref{prop_public}). 
These results follow the intuition that the designer would like to reveal as much information as possible when the payoffs of players are aligned with the system-level objective \cite{Bergemann2019}. {\bf 2)} When the objective is to maximize the agreement between players' actions, we show that no information disclosure is optimal for any LQG game (Proposition \ref{prop_min_dev}). That is, by hiding information players' actions are closer to each other. 
{\bf 3)} If the information designer aims to maximize social welfare and agreement, we identify a critical weight on the agreement term of the objective based on game payoffs below which full information disclosure is preferred to no information disclosure (Propositions \ref{prop:mult_comm} and \ref{prop:mult2}). That is, the benefit of revealing information outweighs the increase in disagreement. {\bf 4)} Numerical solutions to the SDP formulation reveal optimal private signal distributions that outperform both full and no information disclosure schemes.  These contributions mentioned above build on the SDP formulation of the information design problem that considers generic quadratic design objectives in LQG games \cite{Ui2020}, but are distinct in that they provide specific insights about the practically-relevant social welfare and agreement design objectives. 

\noindent {\bf Other related literature:} Other intervention mechanisms, besides information design, include providing financial incentives in the form of taxes and rewards \cite{brown2017studies}, system utility design \cite{Li_Marden}, and nudging or player control during learning dynamics \cite{balcan2012minimally,guers2013informational,riehl2018incentive,das2022average}. 
In contrast to these approaches, the information design framework aims to manage the uncertainties of players so that their expected payoffs align with the objective of a system designer. That is, the system designer does not control agents directly, rather it determines the  information revealed to the players, so that players' evaluation of their payoffs lead to better outcomes from the system designer's perspective. In this sense, there is a limit to the system designer's capability to achieve its goal. This limit determines the value of information.

\vspace{-6pt}
\subsection{Notation}
The $i$th row and $j$th column of matrix $A$ is denoted with $A_{i,j}$.
We use brackets $[A]_{i,j}$ to indicate the $i,j$th submatrix of $A$. For matrices $A\in \reals^{m\times m}$ and $B\in \reals^{m\times m}$, we use $\circ$ to represent the Hadamard product, e.g., $(A \circ B)_{i,j} = A_{i,j} B_{i,j}$.
We use $\bullet$ to represent the Frobenius product, e.g., $A\bullet B = \sum_{i=1}^{m} \sum_{j=1}^{m} A_{i,j}B_{i,j}$. We use $P^m$ and $P^{m}_{+}$ to represent the set of $m \times m$ symmetric and symmetric positive semi-definite matrices, respectively. Trace of a matrix is denoted with $\tr(\cdot)$. $I$ indicates an identity matrix. $\bbone$ is a column vector of all ones.

\section{Information Design in Linear-Quadratic-Gaussian (LQG) Games} \label{sec_model}


\begin{figure}[t]
\centering
\begin{tikzpicture}[level distance=2cm,
  level 1/.style={sibling distance=2.5cm},
  level 2/.style={sibling distance=4cm}]

  \node [cd] {Information Designer ($\zeta$)}
  	child {node [ad] {Player $1$ $(a_{1}, \gamma_{1})$} 	
     edge from parent
       node [below] {$\omega_{1} $}
   }
     child {node {$\ldots$} 	 
   }
       child {node [ad] {Player $n$ $(a_{n},\gamma_{n})$} 
       edge from parent
       node [below] {$\omega_{n}$}
    };
\end{tikzpicture}
\caption{An information designer sends a signal $\omega_{i}$ drawn from information structure $\zeta(\omega| \gamma)$ to each player $i$ who takes action $a_{i}$ in a game with other players under payoff state $\gamma_{i}$.}
\label{inf_design_diagram}
\vspace{-12pt}
\end{figure}

A non-cooperative incomplete information game involves a set of $n$ players belonging to the set $N$, each of which selects actions $a_i \in A_i$ to maximize the expectation of its payoff function $u_i(a,\gamma)$ where $a \equiv (a_{i})_{i \in N} \in A $ and $\gamma \equiv (\gamma_{i})_{i \in N} \in \Gamma$ correspond to an action profile and an unknown payoff state, respectively. Players form expectations about their payoffs based on their signals/types $\omega_i$ about the state given a common prior $\psi$. We represent the incomplete information game by the tuple $G:= \{N, A, \{u_i\}_{i\in N}, \{\omega_i\}_{i\in N}\}$. 

A strategy of player $i$ maps each possible value of the private signal  $\omega_{i} \in \Omega_{i}$ to an action $s_{i}(\omega_{i}) \in A_{i}$, i.e., $s_{i}: \Omega_{i} \rightarrow A_{i}$.  A strategy profile $s = (s_{i})_{i \in N}$ is a BNE with information structure $\zeta$, if it satisfies 
\begin{equation}\label{eq_bne}
E_{\zeta}[u_{i}(s_{i}(\omega_{i}), s_{-i},\gamma )|\omega_{i}] \geq E_{\zeta}[u_{i}(a_{i}', s_{-i},\gamma )|\omega_{i}], 
\end{equation}for all $a_{i}' \in A_{i}, \omega_{i} \in \Omega_{i}, i \in N$ where $s_{-i}=  (s_{j}(\omega_{j}))_{j\neq i}$ is the equilibrium strategy of all the players except $i$, and $E_{\zeta}$ is the expectation operator with respect to the signal distribution $\zeta$ and the prior on the payoff state $\psi$. The above definition ensures that no player has a unilateral profitable deviation from a BNE strategy to another action at any signal realization given the information structure $\zeta$.  

An information designer aims to optimize the expected value of a design objective $f(a,\gamma)$, e.g., social welfare, by deciding on an information structure $\zeta$ from a set of signal generating distributions $Z$, i.e., 
\begin{equation}\label{eq:orig-obj}
\max_{\zeta \in Z} E_{\zeta}[f(s, \gamma)]
\end{equation}
where $s$ is a BNE strategy profile for the game $G$ under the information structure of the game $\zeta$. The information structure of the game  $\zeta(\omega|\gamma)$ is the conditional distribution of $\omega\equiv (\omega_{i})_{i \in N}$ given $\gamma$.
Next, we introduce the two information design objectives that we focus in this paper. 
\begin{example}[Social Welfare]
Social welfare is the sum of individual utility functions, 
\begin{align}\label{eq_soc_welfare}
f(a, \gamma) &= \sum_{i=1}^{n} u_{i}(a,\gamma).
\end{align}
Social welfare is a common design objective used in congestion  \cite{brown2017studies,wu2021value}, global \cite{morris2002social}  or public goods games \cite{alizamir2020healthcare}.

\end{example}
\begin{example}[Agreement]
The information designer would like players to agree by minimizing the deviation of players' actions from the mean action, i.e., by maximizing 
\begin{equation}\label{eq:ssd_obj}
f(a,\gamma) = -\sum_{i=1}^{n}(a_{i} - \bar a)^{2} \text{, where } \bar a=\frac{1}{n}\sum_{i=1}^{n}a_{i},
\end{equation}
where we assume $a_i \in A_i \equiv \reals$. 
The objective is suitable in settings where consensus is desirable but not exactly attainable. 
For instance, this objective can be used in reducing consumption variability in demand response \cite{eksin2015demand}, or coordinated autonomous movement \cite{eksin2014bayesian}. 
\end{example}

Information design follows the given timeline (Fig. \ref{inf_design_diagram}). 
\begin{enumerate}
\item Designer selects $\zeta \in Z$ and notifies all players.
\item Payoff state $\gamma$ is realized.
\item Players observe signals $\{\omega_{i}\}_{i \in N}$ drawn from $\zeta(\omega|\gamma)$. 
\item Players act according to BNE under $\zeta.$
\end{enumerate}

The information designer's problem in \eqref{eq:orig-obj} is intractable for general incomplete information games with continuous actions because it is a linear program with an infinite number of variables \cite{Ui2020}. In this paper, we focus on LQG games that admit a tractable SDP formulation for \eqref{eq:orig-obj} when $f(\cdot)$ is quadratic and signals come from a Gaussian distribution.

\subsection{Linear-Quadratic-Gaussian Games}
In a LQG game, player $i$'s payoff function is quadratic, 
\begin{equation}\label{utility}
u_{i}(a , \gamma ) = - H_{i,i}a_{i}^{2} - 2 \sum_{j \neq i}H_{i,j}a_{i}a_{j} + 2\gamma_{i}a_{i} +d_{i}(a_{-i}, \gamma),
\end{equation}
where $H_{i,j}$ for $i\in N\,,\; j\in N$ are real-valued coefficients with $H_{i,i}>0$, $d_{i}(a_{-i},\gamma)$ is an arbitrary function of the opponents' actions $a_{-i}\equiv \{a_{j}\}_{j \neq i}$ and state $\gamma$, and we have $a \in A \equiv \mathbb{R}^{n} $, and $\gamma \in \Gamma \equiv \mathbb{R}^{n}$. We collect the payoff function coefficients in a matrix $H = [H_{i,j}]\in \reals^{n\times n}$. We note that the function is quadratic in player $i$'s action but it need not be quadratic in others' actions and payoff state as per the term $d_{i}(a_{-i},\gamma)$. Indeed, this term cannot be controlled by player $i$, i.e.,  it does not affect its strategy. Here, we focus on scalar actions, i.e., $a_i \in \reals$. 
\begin{remark}The results in the paper can be extended to cover the case where $a_{i} \in  \mathbb{R}^{m_i}$ for $m_i\in \naturals$, as long as $u_{i}(a, \gamma)$ remains a quadratic in actions. 
\end{remark}

Payoff state $\gamma$ follows a normal distribution $\psi(\mu , \Sigma)$ with mean $\mu\in \mathbb{R}^n$ and covariance matrix $\Sigma$. 
Player $i$ receives a private signal $\omega_{i} \in  \mathbb{R}$. We assume the joint distribution over the random variables $(\omega,\gamma)$ is normal; thus, $\zeta$ is assumed to be a normal distribution. 
Next, we provide two canonical examples of quadratic payoffs. 

\subsubsection{Cournot competition}
Firms determine the production quantities for their goods ($a_{i}$) facing a marginal cost of production ($\gamma_i$) \cite{vives2010information}. The price is a function of the production quantities,  $p_{i}(a)=\vartheta-\varpi a_{i} - \varrho\sum_{j\neq i} a_{j}$ with positive constants $\vartheta$, $\varpi$ and $\varrho$.
The payoff function of the firm $i$ is its profit given by its revenue $a_{i} p_{i}(a) $ minus the cost of production $\gamma_{i}a_{i}$,
\begin{align}\label{eq_cournot}
 u_{i}(a, \gamma) = a_{i} p_{i}(a)-\gamma_{i}a_{i}.  
\end{align}
%

\subsubsection{Beauty Contest Game}
Payoff function of player $i$ is given by
\begin{equation}\label{eq:payoff_beauty}
u_{i}(a, \gamma) = -(1-\beta)(a_{i}-\gamma)^{2} - \beta(a_{i}-\bar{a}_{-i})^{2},
\end{equation}where $\beta\in [0,1]$  and $\bar{a}_{-i}= \sum_{j\neq i} a_{j}/(n-1)$ represents the average action of other players. The first term in \eqref{eq:payoff_beauty} denotes the players' urge for taking actions close to the payoff state $\gamma$. The second term accounts for players' tendency towards taking actions in compliance with the rest of the population. The constant $\beta$ gauges the importance between the two terms. The payoff captures settings where the valuation of a good, e.g., stock, depends not just on the performance of the company but also on what other players think about its value \cite{morris2002social}. 

\subsection{Preliminaries: A SDP Formulation of Information Design Problem given Quadratic Design Objectives}

In this section, we provide preliminary results on the information design problem in LQG games. The first result represents the problem in \eqref{eq:orig-obj} as a SDP with
the decision variable 
$
X:=\begin{bmatrix}
var(a) & cov(a, \gamma)\\
cov(\gamma, a) & var(\gamma)
\end{bmatrix}.
$

\begin{proposition}[\cite{Ui2020}] \label{prop_SDP}
If the objective function $f(a,\gamma)$ is quadratic in its arguments, and the payoff matrix $H$ is such that $H+H^T$ is positive definite, then the information design problem in \eqref{eq:orig-obj} can be restated as the following SDP,
\begin{align}
\max_{X \in P^{2n}_{+}}&  \hspace{-3pt} F\hspace{-1pt} \bullet\hspace{-1pt} X  \hspace{-4pt} = \hspace{-4pt}\max_{X \in P_{+}^{2n}} \hspace{-2pt} \begin{bmatrix}
[F]_{1,1} & [F]_{1,2}\\ [F]_{1,2} & [F]_{2,2}
\end{bmatrix}\hspace{-1pt} \bullet X
\label{eq:mod1}\\
\text{s.t. }& \;R_{k}\bullet X = 0 \quad \forall \; k \in {\{1, . . , n\}},\label{eq:modend} \\
& M_{k,l}\bullet X = \cov(\gamma_{k},\gamma_{l}), \quad\forall\; k,l \in N \text{ with } k\leq l \label{eq:mod2}
\end{align}
where $[F]_{i,j}$ indicates the $n \times n$ block matrix for $i,j\in\{1,2\}$, 
$R_{k} \in P^{2n}$ and $M_{k,l} \in P^{2n}$ are defined as 
\begin{equation*}
[R_{k}]_{i,j}=\begin{cases}
H_{k,k}&if \quad  i = j = k,  \\
H_{k,j}/2  & if \quad  i = k, 1\leq j \leq n, j\neq k, \\
-1/2 &if \quad i = k, j= n + k,\\
H_{k,i}/2 &if\quad  j = k, 1\leq i \leq n, i\neq k \\
-1/2 &if \quad j = k, i= n + k,\\
0 & \text{otherwise,}
\end{cases}
\end{equation*}
and
\begin{equation*}
[M_{k,l}]_{i,j}=\begin{cases}
1/2 \quad\text{ if } k<l, i=n+k, j=n+l\\
1/2 \quad\text{ if } k<l, i=n+l, j=n+k\\
1 \quad\text{ if } k=l, i=n+k, j=n+l\\
0 \quad \text{otherwise.}
\end{cases} 
\end{equation*}
\end{proposition}
This result, due to \cite{Ui2020}, represents the original information design problem \eqref{eq:orig-obj} as the maximization of a linear function of a positive semi-definite matrix $X$ subject to linear constraints. The result leverages the fact there is a unique BNE that is a linear function of the signals whose coefficients can be obtained by solving a set of linear equations in a LQG game with payoff matrix $H$ where $H+H^T\in P^n_+$ \cite{radner1962}. The linear strategies allow a mapping from strategies to signals, which then means selecting the best distribution over the signals is equivalent to selecting the best distribution over the actions subject to the BNE constraints. Accordingly, the selection of the information structure in \eqref{eq:orig-obj} reduces to determining the covariance between the realized actions and payoff states in \eqref{eq:mod1}. Note that we can assume $[F]_{2,2}$ is a zero matrix $\textit{O}_{n\times n}$, because $var(\gamma)$ is given by nature, and cannot be altered by choosing an information structure. Again by leveraging the linear mapping of strategies from signal space to action space, one can express the BNE equations with the set of linear constraints in \eqref{eq:modend}. The set of constraints in \eqref{eq:mod2} assigns the given covariance matrix of the payoff states to the corresponding sub-matrix in $X$, i.e., it is equivalent to $[X]_{2,2} = \var(\gamma)$.  We note that we assume the conditions in Proposition \ref{prop_SDP} hold throughout the paper. 
%

Next, we consider an important special case. 


\begin{definition}[Public Information Structure]\label{defn_public}
A public information structure has $\omega_{1}=....= \omega_{n}$ with probability one. 
The set of public information structures is a subset of the general information structures.
\end{definition}
In the public information design problem, all players receive the same signal, and it is common knowledge that they will receive the same signal.
\noindent We define two important feasible solutions to \eqref{eq:mod1} - \eqref{eq:mod2} (no and full information disclosure).
\begin{definition}[No information disclosure]\label{def_no_info}
No information disclosure refers to the case when there is no informative signal sent to the players. In this case, the equilibrium action profile is given by $a=H^{-1}\mu$. The induced decision variable and the objective value is respectively given by  
\begin{equation}\label{eq_no_info}
X=
\begin{bmatrix}
\textit{O} & \textit{O} \\\textit{O}&\var(\gamma) 
\end{bmatrix}\text{ and } F\bullet X =0.
\end{equation}
\end{definition}

\begin{definition}[Full information disclosure] \label{def_full_info}
The signals sent to the players reveal all elements of payoff state $\gamma$ under full information disclosure. Equilibrium action profile is given by $a=H^{-1}\gamma$. The induced decision variable 
\begin{equation}\label{eq_full_info_solution}
X=\begin{bmatrix}H^{-1}\var(\gamma)(H^{-1})^{T} & H^{-1}\var(\gamma) \\
\var(\gamma)(H^{-1})^{T}  & \var(\gamma)
\end{bmatrix}  
\end{equation} and the objective value is $F\bullet X = F_{H}\bullet \var(\gamma)$ where $F_{H}=  (H^{-1})^{T}([F]_{1,1} + [F]_{1,2}H + H^{T}[F]_{2,1})H^{-1}$.
\end{definition}


Next result states the conditions for the optimality of full information disclosure solution when we consider the set of public information structures. 
\begin{proposition}[Proposition 7,\cite{Ui2020}]\label{pro:full-info}
Let $\var(\gamma) ={DD}^{T}$ such that $D$ is an $n \times k$ matrix of rank $k$ where $k$ is the rank of $\var(\gamma)$.  Assume $D^{T}F_{H}D \neq \textit{O}$ is positive semi-definite. Then, full information disclosure is optimal in the set of public information structures, and no information disclosure is not optimal in the set of general information structures.
\end{proposition}

\begin{remark}
The SDP formulation of the information design problem in \eqref{eq:orig-obj} poses the problem as the determination of a distribution over actions not signals. A natural question is: how can the designer use the solution $X$ and $\phi$ instead of the distribution over signals $\zeta$? As per the information design timeline, when $X$ is decided and $\gamma$ is realized, the designer can draw the suggested actions from $\phi(a|\gamma)$ which has a Gaussian distribution. These suggested actions can be used as coordinating signals instead of the private signals $\omega_i$. 
\end{remark}

\subsection{Design objectives}
In this paper, we focus on two specific quadratic design objectives: social welfare \eqref{eq_soc_welfare} and agreement \eqref{eq:ssd_obj}. According to Proposition 1, we can express the information design problem in \eqref{eq:orig-obj} for these objectives as in \eqref{eq:mod1}. The following are the objective coefficients matrices
 \begin{equation} \label{eq_social_coeffs}
 F^W=\begin{bmatrix}
 -H & I \\ I& \textit{O}
 \end{bmatrix}, \textrm{ and}\;  F^C =\begin{bmatrix}
 \frac{1}{n} \bbone \bbone^T- I & \textit{O} \\ \textit{O}& \textit{O}
 \end{bmatrix},
 \end{equation}
corresponding to \eqref{eq_soc_welfare} and \eqref{eq:ssd_obj}, respectively. We obtain $F^W$ by substituting the quadratic payoffs \eqref{utility} in \eqref{eq_soc_welfare}, and taking the expectation. See Lemma \ref{lemma_obj_min} in the appendix for the derivation of $F^C$.

 \section{Social Welfare Maximization} \label{sec_welfare}


Our first result shows that full information disclosure will be preferred to no information disclosure in social welfare maximization. 
\begin{proposition}\label{prop:full-vs-no}
Assume $H$ is symmetric. Then, full information disclosure never performs worse than no information disclosure for maximizing social welfare objective.
\end{proposition}
\begin{proof}
No information disclosure has the objective value $F\bullet X = 0$ regardless of $F$ as per Definition \ref{def_no_info}. Objective value of full information disclosure is $F^W\bullet X = F_{H}^W\bullet \var(\gamma)$ as per Definition \ref{def_full_info}. Given \eqref{eq_social_coeffs}, $F_{H}^{W} = H^{-1}$. We have $F_{H}^{W} = H^{-1}\succ 0$ because eigenvalues of $H^{-1}$ is equal to reciprocals of eigenvalues of $H$ which are positive because $H$ is positive definite by the assumption that $H+H^{T} \succ 0$ and $H$ is symmetric. The result follows from the fact that $F_{H}^{W} \bullet \var(\gamma)\geq 0$ given $\var(\gamma) \succeq 0$.
\end{proof}
The result implies that no information disclosure cannot be an optimal information structure for social welfare maximization given symmetric payoff coefficients, since it cannot do better than full information disclosure. Next, we show that full information disclosure maximizes social welfare for some important special cases. 

\subsection{Common Payoff State}
 We consider a scenario in which the payoff states are identical for all the players.
\begin{proposition}\label{prop_comm}
 Assume $H$ is symmetric and $\gamma_{i}=\gamma_{j}, \, \forall\; i,j \in N $. Then, full information disclosure ($X$ in \eqref{eq_full_info_solution}) is optimal for social welfare maximization.
\end{proposition}
 \begin{proof} 
 The objective function $f(\cdot)$ with coefficients matrix $F^{W}$ in  \eqref{eq_social_coeffs} is such that $F_{i, n+j}^{W}=0$ $\forall\; i,j \in N$ with $i \neq j$. Moreover, we have $F_{{n+i},{n+j}}^{W}=0,\, \forall\; i,j \in N$. Therefore,
 \begin{equation}\label{eq:co}
 F^{W}\bullet X=\sum_{i=1}^{n}\sum_{j=1}^{n}F_{i,j}^{W}\cov(a_{i},a_{j}) + 2\sum_{i=1}^{n}F_{i,n+i}^{W}\cov(a_{i},\gamma_{i}).
 \end{equation}
Using the BNE condition in \eqref{eq:modend}, which is equivalent to 
\begin{equation}\label{eq:bce2}
\sum_{j\in N}H_{i,j}\cov(a_{i},a_{j}) = \cov(a_{i},\gamma_{i}), \; \forall \; i,j\in N,
\end{equation}
for the corresponding terms in \eqref{eq:co}, we obtain 
 \begin{align}\label{eq_full_symm}
 F^{W}\bullet X &=\sum_{i=1}^{n}\sum_{j=1}^{n}(F_{i,j}^{W}+2F_{i,n+i}^{W}H_{i,j})\cov(a_{i},a_{j}).
 \end{align}
We can express $F^W \bullet X = E \bullet \var(a)$ where we define $E:=[F^W]_{1,1} + [F^W]_{2,1} \circ H + [F^W]_{1,2} \circ H^T$ using \eqref{eq_full_symm}.
Substituting $F^W$ \eqref{eq_social_coeffs} in $E$, we get $E=H^T$. Since $H$ is symmetric, we have $E=H$.
We have that if $E=\kappa H$ for some constant $\kappa>0$, then full information disclosure is optimal under common payoff states (Proposition $9$ in \cite{Ui2020}). In our setting, the condition holds with $\kappa=1$.  
  \end{proof}

 
Proposition \ref{prop_comm} establishes that full information disclosure is the optimal information structure among all possible information structures if the payoff state is common and $H$ is symmetric. 
In the following example, we analyze the discrepancy between the optimal objective value obtained by solving the SDP in \eqref{eq:mod1}-\eqref{eq:mod2} and full information disclosure, as we gradually relax the assumptions of Proposition \ref{prop_comm}. In particular, we allow partially correlated payoff states, and asymmetric game coefficients $H$. 

\begin{figure}[ht]\label{figure_corr}
\centering
 \includegraphics[width=0.7\linewidth]{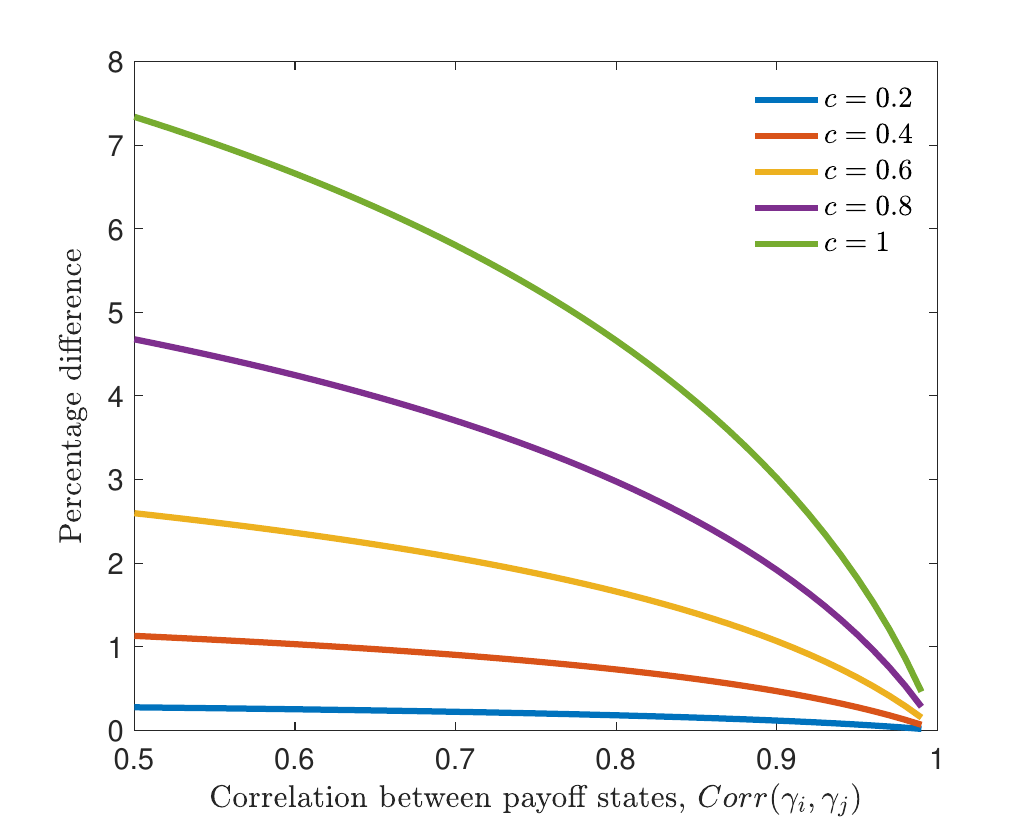}
\caption{\small Percentage difference between optimal objective value \eqref{eq:mod1} and objective value of full information disclosure versus correlation between payoff states.  We consider a game with asymmetric payoffs given by $H_{i,i} = 4$ for $i \in N$, and $H_{i,j} = 1 + c U_{i,j}$ for $i\in N$, and $j \in N \setminus \{i\}$ where $U_{i,j}\in [-1,1]$ is a uniformly distributed random variable` for $i,j \in N$, and $c\in[0,1]$ is a constant determining the magnitude of the asymmetry. The suboptimality of  full  information  disclosure increases  with  growing asymmetry and decreasing correlation.}\label{plot:payoff_asymm}
\end{figure}

\vspace{6pt}
\noindent {\bf Example (Asymmetric payoffs and correlated payoff states):}
Figure \ref{plot:payoff_asymm} shows that full information disclosure becomes increasingly suboptimal as asymmetry grows and correlation between payoff states decreases.
Note that when $Corr(\gamma_i,\gamma_j)=1$, there is a common payoff state and full information disclosure is optimal for symmetric $H$. 
If we consider the beauty contest game with symmetric $H$ and a single stock, full information disclosure on the stock price, i.e, payoff state, is optimal for maximizing social welfare by Proposition \ref{prop_comm}. If we deviate from common payoff state assumption, this means that stock price is not the same for players when they buy the stock. If we deviate from the symmetry assumption, it means the effect of a player $i$'s action on $j$'s payoff is different than the effect of player $j$'s action on $i$'s payoff. In these scenarios, full information disclosure is no longer optimal.


\subsection{Homogeneous LQG games} \label{sec_homogeneous_LQG}



We consider the following payoff matrix $H$:

\begin{equation}\label{H_sym}
H_{i,j}=\begin{cases}
  1\quad\text{ if} \quad  i=j; \; \forall\; i,j \in N \\
 h\quad\text{ if} \quad  i\neq j; \;\forall\; i,j \in N 
\end{cases}
\end{equation}
 in which the off-diagonal terms are identical. {For the Cournot competition \eqref{eq_cournot}, we have a homogeneous payoff matrix with $h= \frac{\varrho}{2\varpi}$ when cost is common, i.e., when $\gamma_i = \gamma_j$ for all $i, j \in N$. For the beauty contest game \eqref{eq:payoff_beauty}, we have an homogeneous payoff matrix with $h=-\frac{\beta}{n-1}$. }

\begin{theorem}\label{thm:sub-super}
Assume $H$ is given in \eqref{H_sym}, and $\tr(\var(\gamma)) \geq 2\sum_{i=1}^{n} \sum_{j\in N \setminus \{i\}} \text{cov}(\gamma_{i},\gamma_{j})$. If $-\frac{1}{n-1}<h<1$, then full information disclosure is optimal for the social welfare maximization objective under general information structures.
 \end{theorem}
 \begin{proof}
 See Appendix \ref{sec_thm_proof} for the proof.
 \end{proof}



\begin{figure}[ht]
\centering
\begin{tabular}{c}
\includegraphics[width=0.7\linewidth]{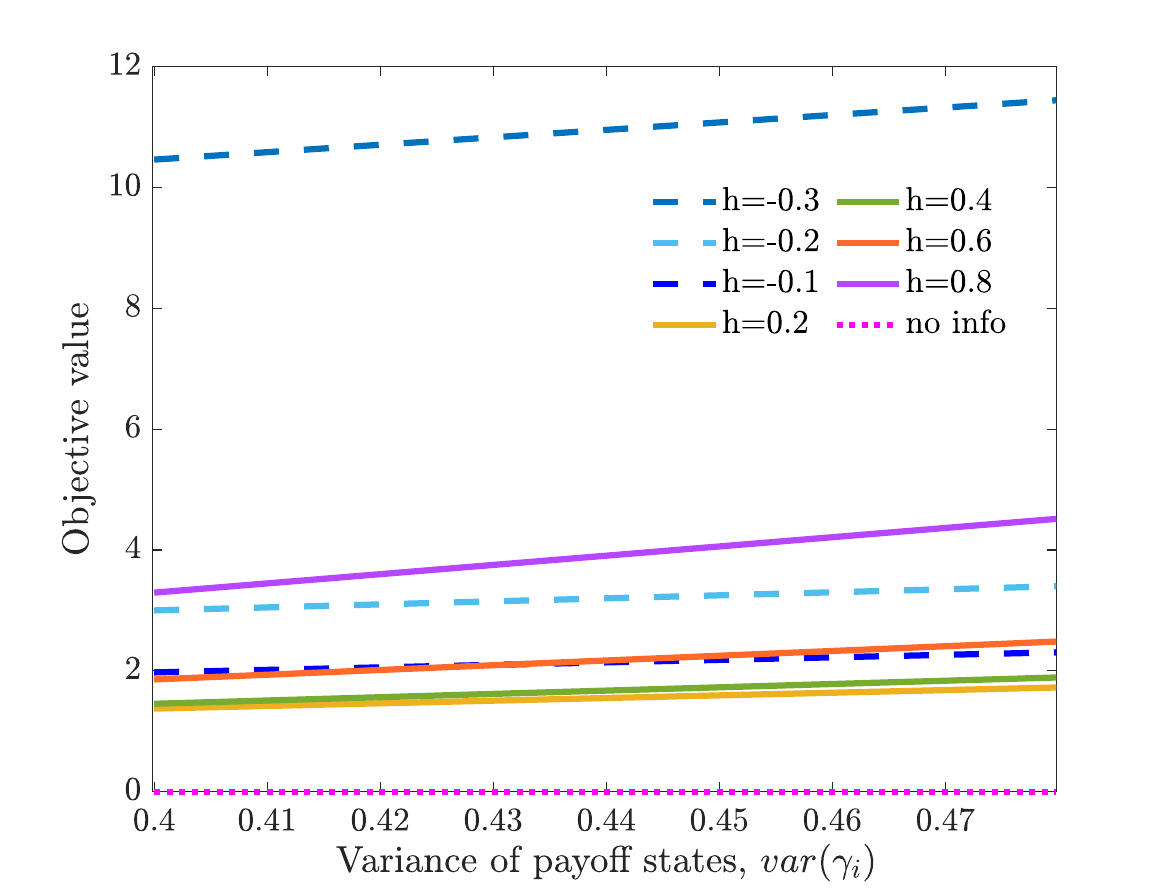}\\
\end{tabular}
\caption{\small Comparison of the social welfare values under full information and no information disclosure. We consider homogeneous games with $\frac{-1}{n-1}<h<1$. We let $\var(\gamma)_{i,j} = 0.2$ for $i\in N$ and $j\in N \setminus \{i\}$ as we vary $\var(\gamma_i)$ for all $i\in N$.    
As $\var(\gamma_i)$ increases, the value of full information disclosure increases compared to no information disclosure.
 }\label{fig:relax}
\end{figure}

Theorem \ref{thm:sub-super} shows that full information disclosure is optimal when the effects of others' actions on payoffs are homogeneous and belong to the given region. We note that the LQG game is submodular if $h> 0$, and it is supermodular if $h<0$. In a submodular game, an increase in a player's action reduces the incentive for other players to increase their actions. In a supermodular game, the effect is reversed, i.e., increasing a player's action increases the incentive for other players to increase their actions---see  \cite{vives2005complementarities} for formal definitions.  Accordingly, social welfare maximization objective is aligned with the incentives of players, when the game is submodular. In contrast, when we have a supermodular game, the optimality of full information disclosure is optimal as long as the effect of another players' actions on a player's action is small, i.e., $h>\frac{-1}{n-1}$. We note that the condition $h>\frac{-1}{n-1}$ stems from the requirement that $H$ needs to be positive definite. Considering the beauty contest game in stock markets with $h=\frac{-\beta}{n-1}<0$, the information disclosure is always optimal because $\beta<1$. In Cournot competition, full information disclosure is optimal as long as $\frac{\varrho}{2\varpi}<1$ according to Theorem \ref{thm:sub-super}.

A sufficient condition for optimality of full information disclosure in Theorem \ref{thm:sub-super} is the diagonal dominance of the covariance matrix of the payoff state. In the following numerical example, we identify that the full information disclosure remains optimal even when the diagonal dominance assumption does not hold in homogeneous LQG games.  

\vspace{6pt}
\noindent {\bf Example (Relaxing the diagonal dominance of $\var(\gamma)$):}
We consider a submodular game among $n=4$ players with homogeneous payoff coefficients with $h$ ranging from $-0.3$ to $0.8$ (see Fig. \ref{fig:relax}).
 When we compare the social welfare value under full information disclosure solution \eqref{eq_full_info_solution} with the optimal solution to the information design problem in \eqref{eq:mod1}-\eqref{eq:mod2}, we find that they are identical for all values of $\var(\gamma_i)\in [0.4, 0.48]$. In this interval of $\var(\gamma_i)$, the diagonal dominance assumption is not satisfied. This suggests that full information disclosure remains optimal even when the diagonal assumption is not satisfied. Fig. \ref{fig:relax} also shows that as the dependence of the payoffs on other players' actions, i.e., $\lvert h\rvert$, increases, objective value for full information disclosure increases. This means the value of revealing information increases as competition increases. 
\subsection{Public information structures} 
Next, we show that full information disclosure maximizes social welfare under public signals.
\begin{proposition}\label{prop_public}
Assume $H$ is symmetric and consider the set of public information structures as the feasible set. Then, full information disclosure maximizes social welfare \eqref{eq_soc_welfare}. 
\end{proposition}
\begin{proof}
From Definition \ref{def_full_info} and $F^W$ in \eqref{eq_social_coeffs}, we have
\begin{equation*}
F_{H}^{W}= (H^{-1})^{T}(-H+ I H + H^{T}I)H^{-1} = H^{-1}.
\end{equation*}

$H^{-1}$ is positive definite because eigenvalues of $H^{-1}$ are equal to reciprocals of eigenvalues of $H$ which are positive. Therefore, $K^{T}F_{H}^W K \neq 0 $ is positive definite for any matrix $K$. The result follows from Proposition \ref{pro:full-info}.
\end{proof}

Together with the previous results in this section, Proposition \ref{prop_public} implies that for full information disclosure to be suboptimal in welfare maximization, the payoff has to include individual payoff states or asymmetric payoff matrix, and the designer has to consider private signals.



 \section{Maximizing Agreement} \label{sec_conformism}

We show that no information disclosure is an optimal information structure that maximizes agreement objective  \eqref{eq:ssd_obj}. 

\begin{proposition}\label{prop_min_dev}
No information disclosure is a maximizer of the objective function in \eqref{eq:ssd_obj} under general information structures.
\end{proposition}
\begin{proof}
The objective coefficients matrix $F^{C}$ has $n-1$ eigenvalues with value $-1$ and $n+1$ eigenvalues with value of $0$. Thus, $F^{C}$ is negative semi-definite. The decision matrix $X$ is positive semi-definite. We deduce that $F^{C}\bullet X \leq 0$. Objective value of no information disclosure is $0$ by \eqref{eq_no_info}; thus, no information disclosure is optimal. 
\end{proof}


Proposition \ref{prop_min_dev} implies that the information designer achieves the maximum similarity between players' actions by revealing uninformative signals to the players. Broadly, hiding information from players is optimal when there is a conflict between the utility functions of the players and the information designer's objective. We compare this with the objective value attained by full information disclosure to provide further intuition for this result. Given $F^C$, we have that $F^C \bullet X = F^C_H \bullet \var(\gamma)$ where $F^C_H = (H^{-1})^T[F^C_{1,1}]H^{-1}$. We know that $\var(\gamma)$ is positive definite, and $F^C_H$ is negative semi-definite because $[F^C]_{1,1}$ is negative semi-definite as per the proof of the Proposition. Thus, we have that full information disclosure can never be better than no information disclosure {for the agreement objective}.

In the context of Cournot competition, we can envision a market regulator that wants to reduce the variability in quantities produced by each firm. The result above states that the designer can achieve minimum variability by not revealing information {about the marginal cost of production}. 

\section{Maximizing Welfare vs. Agreement\label{sec:mult}}
We consider an information design problem in which the designer aims to maximize social welfare and agreement at the same time by considering a weighted combination of \eqref{eq_soc_welfare} and \eqref{eq:ssd_obj}. The objective coefficients matrix is given by 
\begin{equation} \label{mult}
F\hspace{-2pt}:=\hspace{-2pt}((1-\lambda)F^{W} + \lambda F^{C})\hspace{-3pt}=\hspace{-3pt}\begin{bmatrix}
\lambda [F^{C}]_{1,1} - (1-\lambda) H &  (1-\lambda)I \\  (1-\lambda)I& \boldsymbol{O}
\end{bmatrix},
\end{equation}
for weight $\lambda \in  [0,1]$.
%
%
The constant $\lambda$ quantifies the importance of agreement. On one hand full information disclosure is optimal when the design objective is social welfare under common payoff state, homogeneous games, or public signals. On the other hand, no information disclosure is optimal when the objective is to maximize agreement. In the following results we show that full information disclosure remains preferred under public information structures and common payoff states given homogeneous games, if social welfare term gets a large enough weight relative to the agreement term.

%
\begin{proposition}\label{prop:mult_comm}
 Assume $H$ has the form in \eqref{H_sym} with $h\in(0,1)$, and common payoff states $\gamma_{i}=\gamma_{j}, \, \forall\; i,j \in N $. If $\lambda < \frac{1-h}{2-h}$ for $\lambda \in (0,1)$,
full information never performs worse than no information for information design problem with objective coefficients in \eqref{mult}. 
\end{proposition}

\begin{proof}
Following identical steps to Proposition \ref{prop_comm}, we obtain the matrix  $E=[F]_{1,1} + [F]_{2,1} \circ H + [F]_{1,2} \circ H^T$ that provides  $F\bullet X = E \bullet \var(a)$. Substituting in the coefficients from \eqref{mult}, we simplify $E=\lambda [F^C]_{1,1} + (1-\lambda) H$.  
First eigenvalue of $E$ is equal to  $[(n-1)h + 1](1-\lambda)$. The rest of the eigenvalues of $E$ are equal to $-\lambda + (1-\lambda)(1-h)$. $E$ is positive definite because both eigenvalues are greater than zero when $\lambda < \frac{1-h}{2-h}$.
If $E$ is positive definite, then the objective value $E \bullet X_{11} = E\bullet  \var(a)\geq 0$. Thus, full information performs better or the same compared to no information disclosure.
\end{proof}


\begin{proposition}\label{prop:mult2}
Assume $H$ has the form in \eqref{H_sym} with $h\in (0,1)$. If $\lambda < \frac{1-h}{2-h}$ for $\lambda \in (0,1)$, then full information disclosure is  optimal for the information design problem with objective coefficients given in \eqref{mult} under the feasible set of public information structures. 
\end{proposition}
\begin{proof}
We know $F\bullet X= F_H\bullet \var(\gamma)$ where $F_H$ is given in Definition \ref{def_full_info}. Substituting in for the sub-matrices in \eqref{mult}, we have $F_{H} = (H^{-1})^{T} E H^{-1}$,
 where $E= \lambda [F^C]_{1,1} + (1-\lambda) H$ is as in Proposition \ref{prop:mult_comm}. We know from the proof of Proposition \ref{prop:mult_comm} that $E$ is positive definite for $\lambda < \frac{1-h}{2-h}$. Thus, full information disclosure is optimal for public information structures by the fact that $F_H$ is positive semi-definite and by Proposition \ref{pro:full-info}.
\end{proof}

\begin{figure}[ht]
\begin{tabular}{c}
\includegraphics[width=\linewidth]{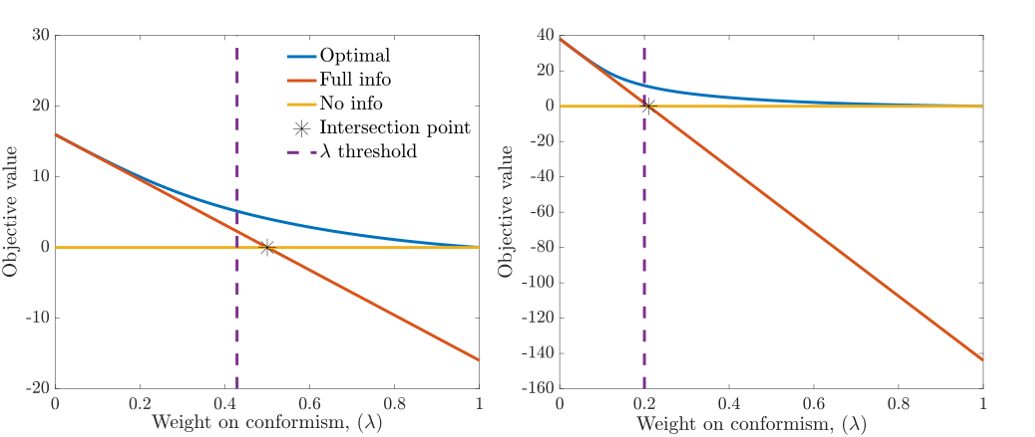}\\
(a) $h=0.25$ \qquad \qquad (b) $h=0.75$
\end{tabular}
\caption{Objective values for optimal, full, no  information disclosure under varying weights $\lambda \in [0,1]$. Optimal information disclosure is obtained by solving \eqref{eq:mod1}-\eqref{eq:mod2} under general information structures. The game payoff coefficients $H$ is as in \eqref{H_sym} with $h \in \{0.25, 0.75\}$. Let $\var(\gamma)$ be such that $\var(\gamma)_{i,i}=4$ for $i\in N$ and $\var(\gamma)_{i,j} =1$ for $i\in N$ and $j\in N\setminus \{i\}$. Full information disclosure is preferred over no information disclosure for larger weight values $\lambda$ than the $\lambda$ threshold given in Proposition \ref{prop:mult2} (dashed line).}
\label{plot:multi} 
\end{figure}



Propositions \ref{prop:mult_comm} and \ref{prop:mult2} specify the threshold for $\lambda$ below which social welfare dominates the agreement term so that no information disclosure can no longer be optimal. It is worth noting that the $\lambda$ threshold for the superiority of full information disclosure are identical in both results. This stems from the fact that we can reduce the objective value evaluation $F\bullet X=E \bullet X$ in both settings. According to the threshold $\lambda < \frac{1-h}{2-h}$, the region in which no information disclosure is not optimal increases to $\lambda \in (0,0.5)$ as $h\to 0^+$. The region in which no information disclosure is not optimal shrinks to $\lambda=0$ as $h\to 1$. 
That is, as the dependence of players' payoffs on others' actions increases, no information disclosure can no longer be ruled out as sub-optimal, unless social welfare maximization is the objective of the designer, i.e., $\lambda = 0$.

Next, we assess the tightness of the threshold for $\lambda$, and the optimality of no and full information disclosures for the class of general information structures in a numerical example.

\vspace{6pt}
\noindent{\bf Numerical example:} Fig. \ref{plot:multi} shows that the region for the weight $\lambda$ where full information information disclosure is preferable by the information designer over no information disclosure under public information structures is larger than the region given by the condition $\lambda< \frac{1-h}{2-h}$. The gap between the analytical threshold (dashed line) and the numerical threshold (shown by $*$) for $\lambda$ decreases as $h$ increases. Fig. \ref{plot:multi} also shows the optimal value achieved by solving the information design problem under general information structures. We observe that general information structures that send partial signals to players perform better than no and full information disclosure for $\lambda \in (0,1)$. 



As mentioned in Section \ref{sec_homogeneous_LQG}, Cournot competition is submodular. Fig. \ref{plot:multi} indicates that as $h= \frac{\varrho}{2\varpi}$ decreases, i.e., the value of information disclosure increases. In other words, in settings where competition is fierce, hiding information is preferred when agreement is a design factor.

\section{Conclusions}
We analyzed information design problem in LQG games given social welfare and agreement as design objectives. We showed that full information disclosure is an optimal solution for welfare maximization when there are common payoff states,  competition among players is homogeneous, or the information designer commits to releasing public signals. When the objective is to induce agreement between players' actions, we showed that hiding information is optimal among general information structures. If the design objective is a weighted combination of social welfare and agreement, there exists a critical weight on agreement term below which full information disclosure is preferred. These results follow the intuition that if the objectives of the information designer and the payoffs of players are in conflict, information designer should blur or hide the information. If the objectives align, the information designer should reveal information. 


\bibliographystyle{ieeetr}
\bibliography{ref}

\appendix
\subsection{Coefficients matrix of the agreement objective}
\begin{lemma}[Agreement Objective]\label{lemma_obj_min} 
The expected value of $f(a,\gamma)$ for \eqref{eq:mod1} can be written as $F^C \bullet X$ where $F^C$ is given in \eqref{eq_social_coeffs}.
\end{lemma}
\begin{proof} Via expanding and regrouping the terms in \eqref{eq:ssd_obj}, 
\begin{align}
E[-\sum_{i=1}^{n}(a_{i} -\bar a)^2] &=  
 \sum_{i=1}^{n} \frac{1-n}{n}E[a_{i}^2] + \frac{2}{n}\sum_{i=1}^{n}\sum_{j=1}^{n}E[a_{i}a_{j}] \label{eq:dev_obj}.
\end{align}
Because $E[a_i]$ is constant for all $i\in N$, we can write \eqref{eq:dev_obj} as $F^C \bullet X$ using the definition of $\var(a)$. 
\end{proof}



\subsection{Proof of Theorem \ref{thm:sub-super}} \label{sec_thm_proof}


We verify that the full information disclosure solution satisfies the KKT conditions.  We denote the dual variables associated with constraints \eqref{eq:modend}, \eqref{eq:mod2} and $X\in P^{2n}_{+}$ by $\overline{\nu} \in \mathbb{R}^{n}$, $\overline{\mu} \in \mathbb{R}^{n(n+1)/2}$ and $\overline{\Gamma}$, respectively. Primal feasibility conditions in \eqref{eq:modend}-\eqref{eq:mod2} are satisfied by full information disclosure.  Next we respectively state the rest of the KKT conditions, i.e., dual feasibility, first order optimality and complementary slackness condition,
  \begin{equation}\label{eq:sub_dual_feas}
 \overline{\Gamma} \in P^{2n}_{+},
 \end{equation}
\begin{equation}\label{first_order}
F^W+\sum_{k=1}^{n}\overline{\nu}_{k} R_{k} + \sum_{k=1}^{n}\sum_{l=1}^{k} \overline{\mu}_{(n-1)k+l}M_{k,l}+ \overline{\Gamma}=0,
\end{equation}
\begin{equation}\label{comp_Sl}
{X}\bullet \overline{\Gamma}=0.
\end{equation}
Let $\overline {X} \in P^{2n}_{+}$ denote the full information disclosure solution as given in \eqref{eq_full_info_solution} to the social welfare maximization problem \eqref{eq_soc_welfare} with coefficients $F^W$. We check whether the above KKT conditions are satisfied by $\overline{X}$. 
We look for a uniform dual variable $\overline{\nu}$, i.e., $\overline \nu_{k} =\nu, \forall \; k \in N$ where $\nu \in \mathbb{R}$, that satisfies \eqref{first_order}. {We} define $\Xi=- \sum_{k=1}^{n}\sum_{l=1}^{k} \overline{\mu}_{(n-1)k+l}M_{k,l}$ in matrix form and assume $\Xi=\mu I, \mu>0$. We can express the dual variable $\overline \Gamma$ using \eqref{first_order} and substituting in \eqref{eq_social_coeffs} for $F^W$,
\begin{equation} \label{eq:comp_2}
    \overline \Gamma = \begin{bmatrix}
(1-\nu) H & (\frac{\nu}{2}-1)I \\ (\frac{\nu}{2}-1)I& \Xi.
\end{bmatrix}
\end{equation}

We use Schur complement to analyze the positive definiteness of $\overline{\Gamma}$ in \eqref{eq:comp_2}.  
A strict version of dual feasibility condition $\overline{\Gamma} \succ 0 $  is satisfied if and only if  $\Xi$ is positive definite and Schur complement 
\begin{equation}\label{eq_schur1}
\overline{\Gamma} / \Xi = (1-\lambda)H - \frac{(\frac{\nu}{2}-1)^{2}I}{\mu}.
\end{equation}
of block matrix $\Xi$  of matrix $\overline{\Gamma} $ is positive definite.
Sum of each row of $\overline{\Gamma} / \Xi$ is equal to $(1-\nu) - (\frac{\nu}{2}-1)^{2}/\mu  + (n-1)(1-\nu)h$. {This} is the first eigenvalue of $\overline{\Gamma} / \Xi$. Rest of the eigenvalues of $\overline{\Gamma} / \Xi$ are equal to $(1-\nu)(1+h) - \frac{({\nu}/{2}-1)^{2}}{\mu}$. We have all of the eigenvalues of $\overline{\Gamma} / \Xi$ positive and $\Xi\succ 0$, when
\begin{equation}\label{eq:gamma-pos}
\mu > \max\{\frac{(\frac{\nu}{2}-1)^{2}}{(1-\nu)(1+h)},0 \}.
\end{equation}
Hence, if $\mu$ satisfies (\ref{eq:gamma-pos}), then $\overline{\Gamma}$  is positive definite.



Next, we show that there exists $\nu\in \reals$ and $\mu$ as in \eqref{eq:gamma-pos} satisfying \eqref{comp_Sl}. We can express the inverse of $H$ in \eqref{H_sym} as follows for $n \geq 3$ 
\begin{equation}\label{eq_H_inverse}
H_{i,j}^{-1}= \begin{cases}
\frac{(n-2)h+1}{-(n-1)h^{2}+(n-2)h + 1 } \;\text{ if} \quad  i=j;\; i,j \in N\\
\frac{-h}{-(n-1)h^{2}+(n-2)h+ 1 } \quad\text{ if} \quad  i\neq j; \; i,j \in N
\end{cases}
\end{equation}
When $X=\overline X$ is given by \eqref{eq_full_info_solution} and $\overline \Gamma$ is as in \eqref{eq:comp_2}, we obtain the following equation by computing the Frobenius product terms within \eqref{comp_Sl} corresponding to each of the four sub-matrices, 
\begin{align}
\overline{X}\bullet& \overline{\Gamma}
= n^{2}(1-\nu)^{2}*(\tau + h\phi) \nonumber
\\&+ 2\frac{(\nu-2)[((2-n)h-1)\tau+h \phi]}{(n-1)h^{2}-(n-2)h - 1} + \mu \tau =0, \label{eq:compslack2}
\end{align}
where we let $\tau=\tr(\var(\gamma))$ and  $\phi=2\sum_{i=1}^{n} \sum_{j\in N \setminus \{i\}} \cov(\gamma_{i},\gamma_{j})$ to simplify the exposition.

Next we show that there exists at least one real root of \eqref{eq:compslack2} $\nu\in \reals$ and $\mu$ as in \eqref{eq:gamma-pos}. If there is a real root, there exists a $\nu \in \mathbb{R}$ satisfying the KKT conditions. 

First, we consider the case $\mu= \frac{(\frac{\nu}{2}-1)^{2}}{(1-\nu)(1+h)} + \epsilon, \epsilon >0$. In this case, \eqref{eq:compslack2} becomes
\begin{align}
\overline{X}\bullet \overline{\Gamma} &=n^{2}(1-\nu)^{2}( \tau + h \phi) +  (\frac{(\frac{\nu}{2}-1)^{2}}{(1-\nu)(1+h)} + \epsilon) \tau \nonumber\\&+ \frac{2(\nu-2)[((2-n)h-1)\tau+h \phi]}{(n-1)h^{2}-(n-2)h - 1} +  =0. \label{eq_comp_slack3}
\end{align}
When we equalize the denominators, \eqref{eq_comp_slack3} becomes a cubic equation in $\nu$. The cubic equation with real coefficients always has at least one real root.  

Secondly, we consider the case $\mu=\epsilon, \epsilon>0$. In this case,  \eqref{eq:compslack2} is a quadratic function of $\nu$
\begin{align}
 a \nu^2 + b \nu + c=0, \label{eq: quad-root}
\end{align}
where we define the constants $a$ , $b$ and $c$ as
\begin{equation}
a= n^{2}(\tau + h\phi) 
\end{equation}
\begin{equation}
b= -2n^{2}(\tau + h\phi) + 2\frac{((2-n)h-1)\tau + h\phi }{(n-1)h^{2}+(2-n)h - 1}
\end{equation}
\begin{equation}
c=  n^{2}(\tau + h\phi) +\frac{-4((2-n)h-1)\tau-4h\phi}{(n-1)h^{2}+(2-n)h - 1} + \epsilon\tau
\end{equation}
We want to show $b^{2}-4ac > 0$, so that there exists a real root. Note that $(n-1)h^{2}-(n-2)h - 1 < 0$ for $\frac{-1}{n-1}<h<1$. 
Also, by our assumption $\tau \geq h\phi $. We can deduce that the discriminant ($b^{2}-4ac$) is positive, i.e.,
\begin{align}
b^{2}&-4ac  
= \frac{8n^{2}(\tau + h\phi) [((2-n)h-1)\tau+h\phi]}{(n-1)h^{2}+(2-n)h - 1} \nonumber \\&+ 4 \left(\frac{((2-n)h-1)\tau + h\phi }{(n-1)h^{2}+(2-n)h - 1}\right)^{2} +n^{2}(\tau + h\phi) \epsilon\tau > 0.
\end{align}
Therefore the roots of \eqref{eq: quad-root} are real. We also need to show at least one of the roots of \eqref{eq: quad-root} ($\nu_{r}$) is such that $\nu_{r} > 1$ so that $\mu = \epsilon$ as per \eqref{eq:gamma-pos}.  We consider the larger root, 
\begin{align}\label{eq: lambda-g-one}
\nu_{r} &=  1 - \frac{((2-n)h-1)\tau + h\phi }{n^{2}(\tau + h\phi) [(n-1)h^{2} +(2-n)h - 1]} \nonumber\\&+\frac{ \sqrt{b^{2}-4ac }}{2a}>1.
\end{align}
We know $a>0$. Also, it can be deduced that the third term in \eqref{eq: lambda-g-one} is greater than the absolute value of the second term in \eqref{eq: lambda-g-one}. Thus, $\nu_{r} > 1$.

\end{document}